\input amstex 
\input amsppt.sty
\magnification 1200

\define\la{\lambda}
\define\p{\Cal P}
\define\Z{\Bbb Z}
\define\C{\Bbb C}
\define\R{\Bbb R}
\define\Ma{\operatorname{Mat}}
\define\x{\frak X}
\define\res{\operatorname{Res}\limits_{u=x}}
\define\r{\operatorname{Res}\limits_{\zeta=x}}
\define\al{\alpha}
\define\be{\beta}
\define\ga{\gamma}
\define\de{\delta} 
\define\wt{\widetilde}
\define\const{\operatorname{const}}
\define\ze{\zeta}

\define\wh{\widehat}

\TagsOnRight
\NoBlackBoxes

\topmatter
\title
Riemann--Hilbert problem and the discrete Bessel kernel
\endtitle
\author Alexei Borodin
\endauthor

\abstract
We use discrete analogs of Riemann--Hilbert problem's methods to derive the discrete Bessel kernel which describes the poissonized Plancherel measures for symmetric groups. To do this we define discrete analogs of a Riemann--Hilbert problem and of an integrable integral operator and show that computing the resolvent of a discrete integrable operator can be reduced to solving a corresponding discrete Riemann--Hilbert problem.

We also give an example, explicitly solvable in terms of classical special functions, when a discrete Riemann--Hilbert problem converges in a certain scaling limit to a conventional one;
the example originates from the representation theory of the infinite symmetric group. 
\endabstract

\endtopmatter

\document
\head 1. Introduction
\endhead
Let $S(n)$ be the symmetric group of degree $n$. It is well known that its
irreducible representations can be naturally parametrized by all the
partitions $\la=(\la_1\ge\la_2\ge\dots\ge \la_k)$ of the number $n$ into
natural summands: $n=\la_1+\la_2+\dots+\la_k$, see, e.g., \cite{JK}.
We shall denote the set of all partitions of $n$ by $\p_n$.
For every partition $\la$, let $\dim\la$ be the dimension of the irreducible 
representation corresponding to $\la$. The well--known Burnside formula for
finite groups implies that
$$
\sum_{\la\in\p_n}\dim^2\la=n!.
$$
Thus, we can construct a probability distribution $M_n$ on $\p_n$ by
setting the weight of $\la\in\p_n$ equal to
$$
M_n(\la)=\frac{\dim^2\la}{n!}.
$$
The $M_n$'s are called the {\it Plancherel distributions}. 
These distributions can be obtained from the uniform distributions on
symmetric groups via the Robinson--Schensted correspondence, see, e.g.,
\cite{Sch}. In particular, the distribution of the largest part $\la_1$ of the random partition of $n$ with respect to the Plancherel distribution
coincides with that of the longest increasing subsequence of random
uniformly distributed permutation of degree $n$, see, e.g., \cite{BDJ1, Appendix}.

Let us organize a new probability distribution $M^\theta$ on the set
$\p=\sqcup_{n>1} \p_n$ of all partitions depending on a positive parameter
$\theta>0$ as follows. Let $\la\in \p_n$. Set
$$
M^\theta(\la)=\frac{\theta^n}{n!}\,e^{-\theta}\cdot M_n(\la).
$$

As was shown in \cite{BOO}, the correlation functions of a point process
naturally attached to $M^\theta$ have determinantal form with a certain kernel
expressed through the $J$ Bessel functions. This powerful fact allowed us
to study
asymptotic properties of the Plancherel distributions and, in particular, to
prove the Baik--Deift--Johansson conjecture \cite{BDJ1}, \cite{BDJ2}, about the asymptotic behavior
of finitely many first parts of the large random partition.\footnote{We refer to \cite{BDJ1} and \cite{AD} for references to earlier works on the subject. Note, in particular the papers \cite{LS}, \cite{VK1}, \cite{VK2}.} The same
correlation kernel has also arisen in \cite{J} in the asymptotics
of orthogonal polynomial ensembles related to the Plancherel measures. See also \cite{BOl3} for a discussion of connections between the approaches of \cite{BOO} and \cite{J}. 

The appearance of Bessel functions in \cite{BOO} seemed rather mysterious. A nice representation theoretic explanation was
suggested in
\cite{Ok} --- for a much wider class of measures on partitions it was shown
that the correlation functions are given by determinantal formulas with a
kernel which has a certain integral representation. In the particular
case of Plancherel distributions this approach leads to the integral representation of
Bessel functions, see \cite{BOk, \S4}.

In this note we will obtain the Bessel functions in another way. We will show that the correlation kernel can be defined as a solution of a certain ``discrete Riemann--Hilbert problem''.

 First, we will define discrete analogs of Riemann--Hilbert problems and integrable integral operators and show how the resolvent of a discrete integrable operator can be obtained from a solution of the corresponding discrete Riemann--Hilbert problem. The result is parallel to the known result in the continuous case, see \cite{IIKS}, \cite{D2}. 
 
As was pointed out by P.~Deift and A.~Its, our setting of the discrete Riemann--Hilbert problems is similar to the pure soliton case in the inverse scattering method, see \cite{BC}, \cite{BDT}, \cite{NMPZ, Ch. III}.

Then, applying discrete analogs of standard
methods of the conventional Rie\-mann--Hilbert problem's techniques, we will obtain a system of differential equations on the matrix elements of the solution of our concrete discrete Riemann--Hilbert problem which will lead to the Bessel equation. 

We will also demonstrate an explicitly solvable example of a discrete Riemann--Hilbert problem (more general then the one that arises from Plancherel measures) and a continuous Riemann--Hilbert problem such that the discrete problem converges to the continuous one in a certain scaling limit. The solutions of these problems are expressed through classical special functions. The example has originated from a problem of harmonic analysis on the infinite symmetric group, see \cite{KOV}, \cite{BOl1}, \cite{BOl2}, \cite{BOl3}.

The paper is organized as follows. \S1 is the introduction. In \S2 we describe a problem whose solution provides the correlation kernel for the poissonized Plancherel measures (the discrete Bessel kernel). In \S3 we review general facts about integrable integral operators and Riemann--Hilbert problems. In \S4 we define discrete analogs of Riemann--Hilbert problems and integrable integral operators and show how computing the resolvent of a discrete integrable operator is reduced to solving a discrete Riemann--Hilbert problem. 
In \S5 we apply the general approach discussed in \S3 to a special case which is relevant for us. In \S6 we do the same in the discrete situation. In \S7 we solve the discrete Riemann--Hilbert problem attached to the problem of \S2 using discrete analogs of standard methods of continuous problems. 
In \S8 we discuss the example of explicitly solvable discrete and continuous Riemann--Hilbert problems mentioned above. 

In the preprint version of this paper there was no general setting for discrete Riemann--Hilbert problems and discrete integrable operators. Only the special case of \S6 was worked out. After the preprint had appeared, Percy Deift suggested a general approach to the discrete situation which we present in \S4. I am very grateful to him for allowing me to reproduce his results in this text and for a number of valuable suggestions.

I would like to thank Grigori Olshanski. Without his constant support and stimulating discussions, this work would never be done. I would also like to thank Alexander Its for explaining the basics of Riemann--Hilbert problems to me and for helpful comments. 

\head 2. Setting of the problem
\endhead

We refer the reader to \cite{BOO} for a detailed exposition of the material of this section.

Let us associate to any Young diagram $\la$ a finite point configuration $X(\la)$ in $\Z'$ as follows. Denote by $(p_1,\dots,p_d\,|\,q_1,\dots,q_d)$ the Frobenius coordinates of $\la$ and set
$$
X(\la)=\left\{-q_1-\frac 12,\dots,-q_d-\frac 12,p_d+\frac 12,\dots,p_1+\frac 12\right\}.
$$

We define the {\it correlation functions} $\rho_k(x_1,\dots,x_k)$, $k=1,2,\dots$, of the measure $M^\theta$ by
$$
\rho_k(x_1,\dots,x_k)=M^\theta\{\la\,|\,\{x_1,\dots,x_k\}\subset X(\la)\},\quad x_1,\dots,x_k\in\Z'.
$$

As was shown in \cite{BOO}, the correlation functions have determinantal form,
$$
\rho_k(x_1,\dots,x_k)=\det [K(x_i,x_j)]_{i,j=1}^k,
$$
with a certain kernel $K(x,y)$ on $\Z'$.  

This kernel can be defined by the formula
$K=L(1+L)^{-1}$, where the kernel $L$ (also on $\Z'$) has a particularly simple form, see \cite{BOO, Proposition 2.3}. Namely, its 
block form corresponding to the partition of $\Z'$ into the set of positive half-integers $\Z'_+$ and 
negative half-integers $\Z'_-$ is
$$
L(x,y)=\left[\matrix
0&\frac {\theta^{\frac {x-y}2}}{\Gamma(x+\frac 12)\Gamma(-y+\frac 12)}
\,\frac 1{x-y}\\
\frac {\theta^{\frac {-x+y}2}}{\Gamma(-x+\frac 12)\Gamma(y+\frac 12)}\,\frac 1{x-y}&0
\endmatrix
\right].
\tag 2.1
$$
The notation means that $L(x,y)=0$ if $x$ and $y$ are of the same sign; 
$$
\gathered
L(x,y)=\frac {\theta^{\frac {x-y}2}}{\Gamma(x+\frac 12)\Gamma(-y+\frac 12)}\,\frac 1{x-y},\quad x\in \Z'_+,\,y\in\Z'_-;\\
L(x,y)=\frac {\theta^{\frac {-x+y}2}}{\Gamma(-x+\frac 12)\Gamma(y+\frac 12)}\,\frac 1{x-y},\quad x\in \Z'_-,\,y\in\Z'_+.
\endgathered
$$

This fact is quite elementary. It is more difficult to obtain an explicit formula for $K$. 
In \cite{BOO} we gave several such formulas. Here is one of them (Theorem 1 of \cite{BOO}).

\proclaim{Theorem 2.1} Define the kernel $K(x,y)$ on  $\Z'=\Z'_+\sqcup\Z'_-$ by the formula
$$
K(x,y)=\sqrt{\theta}\,\left[\matrix 
\frac{J_{x-\frac 12}J_{y+\frac 12}-J_{x+\frac 12}J_{y-\frac 12}}{x-y}&
\frac{J_{x-\frac 12}J_{-y-\frac 12}+J_{x+\frac 12}J_{-y+\frac 12}}{x-y}\\
\frac{J_{-x-\frac 12}J_{y-\frac 12}+J_{-x+\frac 12}J_{y+\frac 12}}{x-y}&
\frac{J_{-x+\frac 12}J_{-y-\frac 12}-J_{-x-\frac 12}J_{-y+\frac 12}}{x-y}
\endmatrix\right]\,,
\tag 2.2
$$
where $J_x=J_x(2\sqrt{\theta})$ is the Bessel function of order $x$ and argument $2\sqrt{\theta}$, and the diagonal entries are determined by the L'Hospital rule. Then
$$
K=\frac{L}{1+L}\,,
$$ 
where $L$ is defined by {\rm(}2.1{\rm)}. 
\endproclaim

We call $K(x,y)$ the {\it discrete Bessel kernel}. 

In \cite{BOO} we gave no explanation why the Bessel functions appear in the picture, there we have just 
verified that the relation $K=L(1+L)^{-1}$ holds. Several such explanations exist by now.
Originally, we obtained the discrete Bessel kernel as a limit of the {\it hypergeometric kernel} 
introduced in \cite{BOl2}, see also \S8. The hypergeometric kernel describes a certain two-parametric family of measures 
on partitions called {\it z--measures} which can be degenerated to the Plancherel measures. K.~Johansson 
independently obtained the restriction of the discrete 
Bessel kernel to $\Z'_+$ as a limit of Christoffel--Darboux kernels for the Charlier polynomials \cite{J}.
A.~Okounkov showed that the Plancherel measures as well as the z--measures fit into a much more general 
infinite--parametric family of measures on partitions, and he gave an integral representation for the 
corresponding correlation kernels \cite{Ok}. In the particular case of the Plancherel measures one 
obtains the integral representation of the Bessel function in this way, see \cite{BOk, \S4, Ex. 1}.

In this note we suggest yet another way to obtain the discrete Bessel kernel. We will consider $L$ and $K$ as discrete analogs of {\it integrable operators} in the sense of \cite{IIKS}, see also \cite{D2}. To my best knowledge, such an object has not been discussed in the literature before.\footnote{Of course, the Christoffel--Darboux kernels for discrete orthogonal polynomials can be considered as discrete integrable operators, but, as far as I know, this observation has never been used so far.}

It is known that if $L$ is an integrable operator in the usual sense then the operator $L(1+L)^{-1}$ is also integrable, and it can be expressed through a solution of a certain {\it Riemann--Hilbert problem} (RHP, for short) associated with $L$, see \cite{IIKS}, \cite{D2}, and \S3. We will define a discrete analog of a Riemann--Hilbert problem and we will show that solving a certain {\it discrete} Riemann--Hilbert problem (DRHP, for short) is equivalent to computing the resolvent of a discrete integrable operator. Furthermore, we will show how to obtain the discrete Bessel kernel by employing discrete analogs of standard methods used in handling conventional (continuous) Riemann--Hilbert problems.

\head 3. Integrable operators and Riemann--Hilbert problems. General approach
\endhead
This section gives a brief review of the formalism of integrable operators and corresponding Riemann--Hilbert problems. 
We shall follow \cite{D2} in our description of the material. 

Let $\Sigma$ be an oriented contour in $\C$. We call an operator $L$ acting in $L^2(\Sigma,|d\zeta|)$ {\it integrable} if its kernel has the form
$$
L(\zeta,\zeta')=\frac {\sum_{j=1}^N f_j(\zeta)g_j(\zeta')}{\zeta-\zeta'},\quad \zeta,\ \zeta'\in\Sigma, 
$$
for some functions $f_j,\ g_j$, $j=1,\dots,N$. We shall always assume that $$\sum_{j=1}^N f_j(\zeta)g_j(\zeta)=0,\quad \zeta\in\Sigma,$$ so that the kernel $L(\zeta,\zeta')$ is nonsingular (this assumption is not necessary for the general theory).

We did not impose any restrictions on the functions $f_i,g_i$ and on the contour $\Sigma$ here, these restrictions depend on a particular problem. For example, one can demand that $f_i,g_i\in L^2(\Sigma,|d\zeta|)\cap L^{\infty}(\Sigma,|d\zeta|)$, and the contour $\Sigma$ is such that the Cauchy Principal Value operator $H$,
$$
(Hh)(\zeta)=\lim_{\varepsilon\to +0}\frac 1{\pi i}\,\int\limits_{\{\zeta'\in\Sigma:|\ze'-\ze|>\varepsilon\}}
\frac{h(\zeta')}{\zeta-\zeta'}\,d\zeta'\,,
$$
is $L^2$--bounded. These restrictions imply, in particular, that 
the operator $L$ is a bounded operator in $L^2(\Sigma,|d\zeta|)$, see \cite{D2}.  

The concept of an integrable operator was first introduced in \cite{IIKS}.

It turns out that for an integrable operator $L$ such that $(1+L)^{-1}$ exists the operator  $K=L(1+L)^{-1}$ is also integrable.
\proclaim{Proposition 3.1 [IIKS]} Let $L$ be an integrable operator as described above and $K=L(1+L)^{-1}$. Then the kernel $K(\zeta,\zeta')$ has the form
$$
K(\zeta,\zeta')=\frac {\sum_{j=1}^N F_j(\zeta)G_j(\zeta')}{\zeta-\zeta'},\quad \zeta,\zeta'\in\Sigma, 
$$
where
$$
F_j=(1+L)^{-1}f_j,\qquad G_j=(1+L^t)^{-1}g_j,\quad j=1,\dots,N.
$$
If $\sum_{j=1}^N f_j(\zeta)g_j(\zeta)=0$ on $\Sigma$, then 
$\sum_{j=1}^N F_j(\zeta)G_j(\zeta)=0$ on $\Sigma$ as well.
\endproclaim

A remarkable fact is that $F_j$ and $G_j$ can be expressed through a solution of an associated Riemann--Hilbert problem. 

As we move along $\Sigma$ in the positive direction, we agree that the $(+)$-side (respectively, $(-)$-side) lies to the left (respectively, right). 

Let $v$ be a map from $\Sigma$ to $\Ma(k,\C)$, $k$ is a fixed integer.

We shall say that a matrix function $m:\C\setminus\Sigma\to \Ma(k,\C)$ is a solution of the RHP $(\Sigma,v)$ if the following conditions are satisfied 
$$
\align
&\bullet\quad  m(\zeta)\text{ is analytic in }\Bbb C\setminus \Sigma,\\
&\bullet\quad m_+(\zeta)=m_-(\zeta)v(\zeta),\  \zeta\in\Sigma,\ \text{where }m_\pm(\zeta)=\lim_{\Sb \zeta'\to \zeta\\ \zeta'\in (\pm)\text{-side}\endSb}m(\zeta'),\\
&\bullet\quad m(\zeta)\to I \text{ as } \zeta\to \infty.
\endalign
$$
The matrix $v(\zeta)$ is called the {\it jump matrix}.

\proclaim{Proposition 3.2 [IIKS]} Let $L$ be an integrable operator as described above and $m(\zeta)$ be a solution of the RHP $(\Sigma, v)$ with
$$
v(\zeta)=I-2\pi i\,f(\zeta)g(\zeta)^t\in \Ma(N,\C),
$$
where
$$
f=\left(f_1,\dots,f_N\right)^t,\qquad
g=\left(g_1,\dots,g_N\right)^t.
$$
Then the kernel of the operator $K=L(1+L)^{-1}$ has the form
$$
K(\zeta,\zeta')=\frac { G^t(\zeta')F(\zeta)}{\zeta-\zeta'},\quad \zeta,\zeta'\in\Sigma, 
$$
where 
$$
F=\left(F_1,\dots,F_N\right)^t,\qquad
G=\left(G_1,\dots,G_N\right)^t
$$
are given by
$$
F(\zeta)=m_\pm(\zeta)\, f(\zeta),\qquad G(\zeta)=(m_\pm ^{t}(\zeta))^{-1}\, g(\zeta).
$$
\endproclaim

\head 4. Discrete integrable operators and discrete Riemann--Hilbert problems. General approach
\endhead
In this section we prove discrete analogs of the two results stated in the previous section. The material of this section is due to Percy Deift. 

Let $\x$ be a discrete locally finite subset of $\C$. We call an operator $L$ acting in $\ell^2(\x)$ {\it integrable} if its matrix has the form
$$
L(x,x')=\cases \dfrac {\sum_{j=1}^N f_j(x)g_j(x')}{x-x'},&\quad x\ne x',\\
0,&\quad x=x',
\endcases
\tag 4.1
$$
for some functions $f_j,\ g_j$ on $\x$, $j=1,\dots,N$, satisfying the relation 
$$
\sum_{j=1}^N f_j(x)g_j(x)=0,\quad x\in\x.
\tag 4.2
$$ 
We will assume that $f_j,g_j\in\ell^2(\x)$ for all $j$. 

We will also call an operator in $\ell^2(\x)$ integrable if its matrix has the form (4.1) for $x\ne x'$ and has arbitrary (not necessarily zero) diagonal values.

Set
$$
f=(f_1,\dots,f_N)^t,\quad g=(g_1,\dots,g_N)^t.
$$
Then (4.2) can be rewritten as $g^t(x)f(x)=0$.
We will also assume that the operator
$$
(Th)(x)=\sum_{x'\in\x,\, x'\ne x}\frac{h(x')}{x-x'}
\tag 4.3
$$
is a bounded operator in $\ell^2(\x)$\footnote{This is an analog of the Cauchy Principal Value operator $H$ introduced in the previous section.}. These restrictions guarantee that $L$ is a bounded operator in $\ell^2(\x)$; a proof of this fact is contained in the proof of Proposition 4.3 below. 

As in the continuous case, it turns out that if the operator $(1+L)$ is invertible then $K=L(1+L)^{-1}=1-(1+L)^{-1}$ is also an integrable operator with possibly nonzero diagonal entries.  
\proclaim{Proposition 4.1}
Let $L$ be an integrable operator as described above and $K=L(1+L)^{-1}$. Then the matrix $K(x,x')$ has the form
$$
K(x,x')=\frac {\sum_{j=1}^N F_j(x)G_j(x')}{x-x'},\quad x\ne x', 
$$
where
$$
F_j=(1+L)^{-1}f_j,\qquad G_j=(1+L^t)^{-1}g_j,\quad j=1,\dots,N,
$$
and
$$
K(x,x)=-\sum_{y\in\x,\,y\ne x}L(x,y)K(y,x).
$$
Moreover, $\sum_{j=1}^N F_j(x)G_j(x)=0$ for any $x\in\x$.
\endproclaim 
\demo{Proof} We will follow the lines of the proof of Proposition 3.1 from \cite{KBI, XIV.1}.
The defining relation $K+LK=L$ reads
$$
K(x,y)+\sum_{t\ne x}\frac {f^t(x)g(t)}{x-t}\,K(t,y)=\frac {f^t(x)g(y)}{x-y}\,.
$$
Assume $x\ne y$. Multiplying both sides by $(x-y)$ and using the relation 
$$
\frac {x-y}{x-t}=1+\frac {t-y}{x-t}
$$ 
we get
$$
(x-y)K(x,y)+\sum_{t\ne x} \frac{f^t(x)g(t)}{x-t}\, (t-y)K(t,y)=f^t(x)\left(g(y)-\sum_{t\ne x}g(t)K(t,y)\right).
$$  
Thanks to (4.2), the restriction $t\ne x$ in the second summation above can be removed, and we obtain
$$
\sum_t\left(\delta(x-t)+L(x,t)\right)(t-y)K(t,y)=
f^t(x)\left((1-K^t)g\right)(y).
$$
Applying $(1+L)^{-1}$ to both sides we see that
$(x-y)K(x,y)=F^t(x)G(y)$, where 
$$
F=(F_1,\dots,F_N)^t,\quad G=(G_1,\dots,G_N)^t.
$$
This proves the first claim of the proposition. 

Since the diagonal entries of $L$ are zeros, the relation $K+LK=L$ implies
$$
K(x,x)=-\sum_{y\ne x} L(x,y)K(y,x)+L(x,x)=-\sum_{y\ne x} L(x,y)K(y,x),
$$ 
and this proves the second claim of the proposition. 

The relation $F^t(x)G(x)=0$ will be proved later, see (4.12).
\qed
\enddemo

\example{Remark 4.2}
It is not difficult to see that for an integrable operator $V$ with arbitrary diagonal entries bounded from $-1$, and such that $(1+V)$ is invertible, the operator $W=V(1+V)^{-1}$ is integrable. Indeed, if we denote by $V_d$ the diagonal part of the matrix $V$ and by $V_0$ the off--diagonal part of $V$ (that is, $V=V_d+V_0$) then it is easily seen that
$$
W=V_d\,(1+V_d)^{-1}+\wt W (1+V_d)^{-1},
$$ 
where 
$$
\wt W=\frac{(1+V_d)^{-1}V_0}{(1+(1+V_d)^{-1}V_0)^{-1}}
$$
is integrable by Proposition 4.1. This implies that $W$ is also integrable. 
\endexample

Similarly to the continuous case, $F_j$ and $G_j$ from Proposition 4.1 can be expressed through a solution of an associated {\it discrete} Riemann--Hilbert problem. 

Let $w$ be a map from $\x$ to $\Ma(k,\C)$, $k$ is a fixed integer.

We shall say that a matrix function $m:\C\setminus\x\to \Ma(k,\C)$ 
with simple poles at the points $x\in\x$ is a solution of the DRHP $(\x,w)$ if the following conditions are satisfied 
$$
\align
&\bullet\quad  m(\zeta)\text{ is analytic in }\Bbb C\setminus \x,\\
&\bullet\quad 
\r m(\ze)=\lim_{\zeta\to x}\left(m(\zeta)w(x)\right),\quad x\in\x,\\
&\bullet\quad m(\zeta)\to I \text{ as } \zeta\to \infty.
\endalign
$$
We will also call $w(x)$ the {\it jump matrix}.

If the set $\x$ is infinite, the last condition must be made more precise. Indeed, a function with poles accumulating at infinity cannot have asymptotics at infinity. One way to precise the condition is to require the uniform asymptotics on a sequence of expanding contours, for example, on a sequence of circles $|\ze|=a_k$, $a_k\to+\infty$. 

In order to guarantee the uniqueness of solutions
of the DRHPs considered below we will always assume that there exists a sequence of expanding contours such that the distance from these contours to the set $\x$ is bounded from zero, and we will require a solution $m(\zeta)$ to have the proper asymptotic behavior on these contours.  

The setting of the DRHP above is very similar to the pure soliton case in the inverse scattering method, see \cite{BC}, \cite{BDT}, \cite{NMPZ, Ch. III}.

\proclaim{Proposition 4.3} Let $L$ be an integrable operator as described above and $m(\zeta)$ be a solution of the DRHP $(\x, w)$ with
$$
w(x)=-f(x)g(x)^t\in \Ma(N,\C).
$$
Then the matrix $K=L(1+L)^{-1}$ has the form
$$
K(x,x')=
\cases\dfrac {G^t(x')F(x)}{x-x'},\quad &x\ne x',\\
 G^t(x)\lim\limits_{\zeta\to x}\left(m'(\zeta)
 \,f(x)\right), &x=x',
\endcases
$$
where $m'(\ze)=\dfrac{dm(\ze)}{d\ze}$, and
$$
F(x)=\lim_{\zeta\to x}\left(m(\zeta)\, f(x)\right),\qquad G(x)=\lim_{\zeta\to x}\left((m^{t}(\zeta))^{-1}\, g(x)\right).
$$
\endproclaim

Before proving this claim let us prove the following useful statements.
\proclaim{Lemma 4.4} If $m(\ze)$ is a solution of the DRHP $(\x,w)$, and $w^2(x)=0$ for some $x\in \x$ then the function
$$
m(\ze)\left(I+\frac{w(x)}{\ze-x}\right)
\tag 4.4
$$
is analytic in a neighborhood of $x$.
\endproclaim
\demo{Proof}
In a neighborhood of $x$ we have ($m(\ze)$ has a simple pole at $x$)
$$
m(\ze)=\frac {A_x}{\ze-x}+B_x+O(\ze-x),
$$
where $A_x$ and $B_x$ are constant matrices. 
Then
$$\r m(\ze)=A_x, \quad m(\ze)w(x)=\frac{A_xw(x)}{\ze-x}+B_xw(x)+O(\ze-x).
$$ 
The residue condition of the DRHP implies that
$$
A_xw(x)=0, \quad A_x=-\lim_{\ze\to x}m(\ze)w(x)=-B_xw(x).
$$
Hence,
$$
m(\ze)=B_x\left(-\frac{w(x)}{\ze-x}+I\right)+O(\ze-x)
$$
and, recalling that $w^2(x)=0$, we get
$$
m(\ze)\left(I+\frac{w(x)}{\ze-x}\right)=B_x+O(1)
$$
which is analytic near $x$.\qed
\enddemo

Note that the jump matrix $w(x)=-f(x)g(x)^t$ satisfies the condition $w^2(x)=0$ at every point $x\in\x$ because of (4.2).

\proclaim{Lemma 4.5}
Let $m_1(\ze)$ and $m_2(\ze)$ be solutions of DRHPs with the same jump matrix $w(x)$, $w^2(x)=0$ for all $x\in\x$, and possibly different asymptotic conditions at infinity. If $m_2(\ze)$ is invertible for $\ze\in\C\setminus \x$ then $m_1(\ze)(m_2(\ze))^{-1}$ is analytic in $\C$. 
\endproclaim
\demo{Proof}
Since $w^2(x)=0$, the determinant of the matrix $\left(I+\frac{w(x)}{\ze-x}\right)$ is identically equal to 1.
In particular, this matrix is invertible.
 
If $m_1(\ze)$ and $m_2(\ze)$ are solutions of DRHPs with the same jump matrix $w(x)$ and $m_2(\ze)$ is invertible then 
$$
\gathered
m_1(\ze)(m_2(\ze))^{-1}=m_1(\ze)\left(I+\frac{w(x)}{\ze-x}\right)
\left(I+\frac{w(x)}{\ze-x}\right)^{-1}(m_2(\ze))^{-1}\\=
\left(m_1(\ze)\left(I+\frac{w(x)}{\ze-x}\right)\right)
\left(m_2(\ze)\left(I+\frac{w(x)}{\ze-x}\right)\right)^{-1}
\endgathered
$$
is analytic near $x\in\x$ by Lemma 4.4. Since $x\in\x$ is arbitrary, $m_1(\ze)(m_2(\ze))^{-1}$ is analytic in $\C$. \qed
\enddemo

\proclaim{Corollary 4.6} A DRHP with the jump matrix of the form $w(x)$, $w^2(x)=0$ for all $x\in\x$, and arbitrary invertible asymptotics at infinity whose any solution is invertible in $\C\setminus\x$ has at most one solution.
\endproclaim
\demo{Proof}
Indeed, if $m_1(\ze)$ and $m_2(\ze)$ are two solutions then
$m_1(\ze)(m_2(\ze))^{-1}$ has no singularities by Lemma 4.5 and is asymptotically equal to $I$ at infinity. By the Liouville theorem, $m_1(\ze)\equiv m_2(\ze)$.\qed
\enddemo
\proclaim{Lemma 4.7} If $m(\ze)$ is a solution of a DRHP with the jump matrix $w(x)$, $w^2(x)=0$ for all $x\in\x$, and asymptotics $I$ at infinity then $\det m(\ze)\equiv 1$. Hence, such DRHP has at most one solution. 
\endproclaim
\demo{Proof} 
Lemma 4.4 implies that
$$
\det m(\ze)=\det \left[m(\ze)\left(I+\frac{w(x)}{\ze-x}\right)\right]
$$
is analytic near $x\in\x$. Since $x\in\x$ is arbitrary, $\det m(\ze)$ is analytic in $\C$, and the asymptotic condition $m(\ze)\to I$ as $\ze\to\infty$, implies that $\det m(\ze)\equiv 1$.\qed
\enddemo  

\demo{Proof of Proposition 4.3}
The proof is similar to that of Proposition 3.1 from \cite{D2}. It is based on the following {\it commutation formula}: if $B_1$ and $B_2$ are Banach spaces then for any bounded operators $D:B_1\to B_2$ and $E:B_2\to B_1$,
$$
\frac{\la}{DE+\la}+D\,\frac 1{ED+\la}\, E=1
$$
in the sense that if $(-\la)\ne 0$ lies in the resolvent set of $ED$ then $(-\la)$ lies in the resolvent set of $DE$ and $(DE+\la)^{-1}=\la^{-1}(1-D(ED+\la)^{-1}E)$, see e.g. \cite{Sak}.
This simple result turns out to have a large number of applications in mathematical physics, see \cite{D1}. 

Let $R_f$ denote the map of right multiplication by the column $N$--vector $f$ taking row $N$--vector functions to scalar functions:
$$
(R_f h)(x)=h(x)f(x)=\sum_{i=1}^Nh_i(x)f_i(x),\quad h=(h_1,\dots,h_N),
$$
and let $R_{g^t}$ denote the map of right multiplication by the row $N$--vector $g^t$ taking scalar functions to row $N$--vector functions:
$$
(R_{g^t}h)(x)=h(x)g^t(x)=(h(x)g_1(x),\dots ,h(x)g_N(x)).
$$ 
By applying $R_f$ to a $N\times N$ matrix valued function we will mean the application of $R_f$ to every row of the matrix and getting a column $N$--vector function as the result. Similarly, by applying $R_{g^t}$ to a column $N$--vector function means the application of $R_{g^t}$ to every coordinate of the vector and getting an $N\times N$ matrix valued function as the result.   

Now observe that the operator $L$ is of the form $L=DE$, where $D=R_f$, $E=TR_{g^t}$, and $T$ was introduced in (4.3). 
The operator $D$ maps the space $B_1$ of row $N$--vector functions on $\x$ with coordinates from $\ell^2(\x)$ to $B_2=\ell^2(\x)$, and $E$ maps $B_2$ to $B_1$. The operator $T$ acts in the space of row $N$--functions coordinatewise.

This implies that $L$ is a bounded operator in $\ell^2(\x)$ if the functions $f_j, g_j$ are bounded and the operator $T$ is bounded. Note that for the boundedness of $L$ instead of the boundedness of $T$ we could require the boundedness of $E=TR_{g^t}$. 

Let us apply the commutation formula for $\la=1$. Since $(1+L)$ is invertible, we have (using Proposition 4.1)
$$
\gathered
F=(1+L)^{-1}f=(1+DE)^{-1}f=f-D\,\frac 1{ED+1}\, Ef
\\=f-R_f\,\frac 1{ED+1}\,\,TR_{g^t}f=f-R_f\,\frac 1{ED+1}\,\,T\,fg^t 
\\=f-R_f\,\frac 1{ED+1}\,ED \,{I}=f-R_f\left(1-\frac 1{ED+1}\right)\,{I}
\\=R_f(ED+1)^{-1}{I}.
\endgathered
$$
Here $I$ denotes the $n\times n$ matrix valued function which is identically equal to 1. 

Set $\mu=(ED+1)^{-1}{I}$. Note that all matrix elements of $(\mu-I)$ are in $\ell^2(\x)$. 
We have just showed that
$$
F=(1+L)^{-1}f=R_f(ED+1)^{-1}{I}=\mu f, 
\tag 4.5
$$
where $\mu={I}-ED\,\mu$, or
$$
\mu(x)={I}-\sum_{y\in\x,\,y\ne x}\frac{\mu(y)(fg^t)(y)}{x-y}.
\tag 4.6
$$
Since 
$$
L^t(x,x')=\cases \dfrac {\sum_{j=1}^N f_j(x')g_j(x)}{x'-x},&\quad x\ne x',\\
0,&\quad x=x',
\endcases
$$
can be obtained from $L(x,x')$ defined by (4.1) by the change $f\to -g$, $g\to f$, we get
$$
G=(1+L^t)^{-1}g=\wt \mu\, g,
\tag 4.7
$$
where
$$
\wt \mu(x)={I}+\sum_{y\in\x,\,y\ne x}\frac{\wt \mu(y)(gf^t)(y)}{x-y}.
\tag 4.8
$$
Set 
$$
m(\ze)=I-\sum_{y\in\x}\frac{\mu(y)(fg^t)(y)}{\ze-y}, \qquad \wt m(\ze)=I+\sum_{y\in\x}\frac{\wt \mu(y)(gf^t)(y)}{\ze-y}.
\tag 4.9
$$
cf. (4.6), (4.8). We have
$$
\r m(\ze)=-\mu(x) (fg^t)(x)=-\lim_{\ze\to x}\left(m(\ze) (fg^t)(x)\right).
$$
Indeed, 
$$
m(\ze)(fg^t)(x)=(fg^t)(x)-\sum_{y\in\x,\, y\ne x}\frac{\mu(y)(fg^t)(y)(fg^t)(x)}{\ze-y}
$$
as $(fg^t)^2(x)=0$ by (4.2).

Thus, $m(\ze)$ is a (unique by Lemma 4.7) solution of the DRHP $(\x,w)$ with $w(x)=-f(x)g(x)^t$ (the asymptotic condition immediately follows from (4.9)). Similarly, using (4.5) we get
$$
F(x)=\mu(x)f(x)=\lim_{\ze\to x}\left(m(\ze)f(x)\right).
\tag 4.10
$$
Same arguments show that $\wt m(\ze)$ is a solution of the DRHP
$(\x,\,g(x)f(x)^t)$, and
$$
G(x)=\wt \mu(x)g(x)=\lim_{\ze\to x}\left(\wt m(\ze)g(x)\right).
\tag 4.11
$$
By Lemma 4.4, the functions
$$
m(\ze)\left(I+\frac{f(x)g(x)^t}{\ze-x}\right) \text{      and       }
\wt m(\ze)\left(I-\frac{g(x)f(x)^t}{\ze-x}\right)
$$
are analytic near $x$. Then, using the relation $(f(x)g(x)^t)^2=0$, we have
$$
\gathered
m(\ze)\wt m^t(\ze)=m(\ze)\left(I+\frac{f(x)g(x)^t}{\ze-x}\right)
\left(I-\frac{f(x)g(x)^t}{\ze-x}\right)\wt m^t(\ze)
\\=\left(m(\ze)\left(I+\frac{f(x)g(x)^t}{\ze-x}\right)\right)
\left(\wt m(\ze)\left(I-\frac{g(x)f(x)^t}{\ze-x}\right)\right)^t
\endgathered
$$ 
which is analytic near $x$. Since $x\in\x$ is arbitrary, $m(\ze)\wt m^t(\ze)$ is analytic in $\C$, and the asymptotic condition of DRHP
implies that $\wt m^t(\ze)=(m(\ze))^{-1}$. Thus, (4.11) turns into
$$
G(x)=\lim_{\ze\to x}\left((m^t(\ze))^{-1}g(x)\right).
$$

Note that by (4.10) and (4.11) we have
$$
G(x)^tF(x)=\lim_{\ze\to x}\left(g(x)^t (m(\ze))^{-1}m(\ze) f(x)\right)=
g(x)^tf(x)=0
\tag 4.12
$$
which proves the last claim of Proposition 4.1.

Finally, on the diagonal by Proposition 4.1, (4.2), (4.5), (4.9) we get 
$$
\gathered
K(x,x)=-\sum_{y\in\x,\, y\ne x}L(x,y)K(y,x)
=-\sum_{y\in\x,\, y\ne x}\frac {g(y)^tf(x)}{x-y}\,\frac{G(x)^tF(y)}{y-x}
\\=G(x)^t\left(\sum_{y\in\x,\, y\ne x}\frac{F(y)g(y)^t}{(x-y)^2}\right)f(x)=G(x)^t\left(\sum_{y\in\x,\, y\ne x}\frac{\mu(y)f(y)g(y)^t}{(x-y)^2}\right)f(x)\\=
G(x)^t\lim_{\ze\to x}\left(m'(\ze) f(x)\right). \qed
\endgathered
$$
\enddemo

\head
5. Integrable operators and Riemann--Hilbert problems. Special case
\endhead
In this section we apply the general formalism of \S3 to a much more special situation.

Let $\Sigma=\Sigma_I\cup\Sigma_{II}$ be a union of two contours, and assume that the operator $L$ in the block form corresponding to this splitting is as follows
$$
L(x,y)=\left[\matrix 0&\frac {h_{I}(x)h_{II}(y)}{x-y}\\
\frac {h_{I}(y)h_{II}(x)}{x-y}&0\endmatrix\right]
\tag 5.1
$$
for some functions $h_I(\,\cdot\,)$ and $h_{II}(\,\cdot\,)$ defined on $\Sigma_I$ and $\Sigma_{II}$, respectively.

Then the operator $L$ is integrable with $N=2$. Indeed,
$$
L(x,y)=\frac{f_1(x)g_1(y)+f_2(x)g_2(y)}{x-y},\qquad x,y\in \Sigma,
$$
where
$$
f_1(x)=\cases h_I(x),&x\in\Sigma_I,\\0,&x\in\Sigma_{II},\endcases\quad
f_2(x)=\cases 0,&x\in\Sigma_I,\\h_{II}(x),&x\in\Sigma_{II},\endcases
$$$$
g_1(x)=\cases 0,&x\in\Sigma_I,\\h_{II}(x),&x\in\Sigma_{II}, \endcases\quad
g_2(x)=\cases h_I(x),&x\in\Sigma_I,\\0,&x\in\Sigma_{II}.\endcases
$$
Then the jump matrix $v(x)$ of the corresponding RHP has the form
$$
v(x)=\cases \left(\matrix1&-2\pi i\,h_I^2(x)\\0&1\endmatrix\right),&x\in\Sigma_I,\\
\left(\matrix 1&0\\-2\pi i\, h_{II}^2(x)&1\endmatrix\right),&x\in\Sigma_{II}.
\endcases
\tag 5.2
$$
It can be easily seen that the RHP in such a situation is equivalent to the following set of conditions:

\noindent$\bullet$  matrix elements $m_{11}$ and $m_{21}$ are holomorphic in $\C\setminus\Sigma_{II}$; 

\noindent$\bullet$  matrix elements $m_{12}$ and $m_{22}$ are holomorphic in $\C\setminus\Sigma_{I}$; 

\noindent$\bullet$  on $\Sigma_{II}$ the following relations hold
$$
\gathered
{m_{11}}_+(x)-{m_{11}}_-(x)=-2\pi i\,h_I^2(x)m_{12}(x),\\
{m_{21}}_+(x)-{m_{21}}_-(x)=-2\pi i\,h_I^2(x)m_{22}(x);
\endgathered
$$
\noindent$\bullet$  on $\Sigma_{I}$ the following relations hold
$$
\gathered
{m_{12}}_+(x)-{m_{12}}_-(x)=-2\pi i\,h_{II}^2(x)m_{11}(x),\\
{m_{22}}_+(x)-{m_{22}}_-(x)=-2\pi i\,h_{II}^2(x)m_{21}(x);
\endgathered
$$

\noindent$\bullet$ $m(x)\sim 1$ as $x\to\infty$.

According to Proposition 3.2, the kernel $K(x,y)$ in the block form corresponding to the splitting $\Sigma=\Sigma_I\cup\Sigma_{II}$ looks as follows
$$
\multline
K(x,y)=\\ \left[\matrix
\frac{h_I(x)h_I(y)(-m_{11}(x)m_{21}(y)+m_{21}(x)m_{11}(y))}
{x-y}&
\frac{h_I(x)h_{II}(y)(m_{11}(x)m_{22}(y)-m_{21}(x)m_{12}(y))}
{x-y}\\
\frac{h_{II}(x)h_I(y)(m_{22}(x)m_{11}(y)-m_{12}(x)m_{21}(y))}
{x-y}&
\frac{h_{II}(x)h_{II}(y)(-m_{22}(x)m_{12}(y)+m_{12}(x)m_{22}(y))}
{x-y}
\endmatrix\right].
\endmultline
$$

\head 6. Discrete integrable operators and discrete Riemann--Hilbert problems. Special case
\endhead

Similarly to \S5, we apply the general approach of \S4 to a special case.

Let $\x$ be a discrete locally finite subset of $\C$ and let $\x=\x_I\sqcup\x_{II}$ be its splitting into two disjoint parts. Let $h_{I}(\,\cdot\,)$, $h_{II}(\,\cdot\,)$ be two functions defined on $\x_I$ and $\x_{II}$, respectively. We will assume that
$h_{I}\in\ell^2(\x_I)$, $h_{II}\in\ell^2(\x_{II})$. 

Consider a matrix $L$ of size $\x\times \x$ which in the block form corresponding to the splitting $\x=\x_I\sqcup\x_{II}$ looks as follows, cf. (5.1),
$$
L(x,y)=\left[\matrix 0&\frac {h_{I}(x)h_{II}(y)}{x-y}\\
\frac {h_{I}(y)h_{II}(x)}{x-y}&0\endmatrix\right].
$$
This matrix defines a bounded operator in $\ell^2(\x)$ if, for example, the operator $T$ defined in (4.3) is $\ell^2$--bounded. 

Let us assume that the operator $(1+L)$ is invertible. This is automatically true if $h_{I}$ and $h_{II}$ are real valued, then $L^*=-L$, and $-1$ cannot belong to the spectrum of $L$.

As in \S5, the operator $L$ is integrable with $N=2$. Indeed,
$$
L(x,y)=\frac{f_1(x)g_1(y)+f_2(x)g_2(y)}{x-y},\qquad x,y\in \x,
$$
where
$$
f_1(x)=\cases h_I(x),&x\in\x_I,\\0,&x\in\x_{II},\endcases\quad
f_2(x)=\cases 0,&x\in\x_I,\\h_{II}(x),&x\in\x_{II},\endcases
$$$$
g_1(x)=\cases 0,&x\in\x_I,\\h_{II}(x),&x\in\x_{II}, \endcases\quad
g_2(x)=\cases h_I(x),&x\in\x_I,\\0,&x\in\x_{II}.\endcases
$$
The jump matrix $w(x)=-f(x)g(x)^t$ of the corresponding DRHP has the form
$$
w(x)=\cases \left(\matrix 0&-h_I^2(x)\\0&0\endmatrix\right),&x\in\x_I,\\
\left(\matrix 0&0\\-h_{II}^2(x)&0\endmatrix\right),&x\in\x_{II}.
\endcases
$$

It is readily seen that the DRHP is equivalent to the following conditions:

\noindent$\bullet$  matrix elements $m_{11}$ and $m_{21}$ are holomorphic in $\C\setminus\x_{II}$; 

\noindent$\bullet$  matrix elements $m_{12}$ and $m_{22}$ are holomorphic in $\C\setminus\x_{I}$; 

\noindent$\bullet$  $m_{11}$ and $m_{21}$ have simple poles at the points of $\x_{II}$, and for $x\in\x_{II}$
$$
\gathered
\res m_{11}(u)=-h_I^2(x)m_{12}(x),\\
\res m_{21}(u)=-h_I^2(x)m_{22}(x);
\endgathered
$$
\noindent$\bullet$  $m_{21}$ and $m_{22}$ have simple poles at the points of $\x_{I}$, and for $x\in\x_{I}$
$$
\gathered
\res m_{12}(u)=-h_{II}^2(x)m_{11}(x),\\
\res m_{22}(u)=-h_{II}^2(x)m_{21}(x);
\endgathered
$$

\noindent$\bullet$ $m(u)\sim 1$ as $u\to\infty$.

Let us indicate how this setting of the problem can be obtained from the continuous case. A continuous Riemann--Hilbert problem is equivalent to a system of integral equations. For the special type of Riemann--Hilbert problems described in \S5 the system takes the form
$$
\alignedat{2}
&\int_{\Sigma_{I}}\frac{h_I^2(x){m_{11}}(x)}{x-y}\,dx={m_{12}}(y),\quad& \int_{\Sigma_{I}}\frac{h_I^2(x){m_{21}}(x)}{x-y}\,dx=1+{m_{22}}(y),\\ 
&\int_{\Sigma_{II}}\frac{h_{II}^2(x){m_{12}}(x)}{x-y}\,dx=1+{m_{11}}(y),\quad &
\int_{\Sigma_{II}}\frac{h_{II}^2(x){m_{22}}(x)}{x-y}\,dx={m_{21}}(y).
\endalignedat
$$

For the first two equations $y$ belongs to $\C\setminus\Sigma_{I}$, and for the last two $y$ belongs to $\C\setminus\Sigma_{II}$.

A natural discrete analog of this system is the following one
$$
\alignedat{2}
&\sum_{x\in\x_{I}}\frac{h_I^2(x){m_{11}}(x)}{x-y}={m_{12}}(y),\quad &
\sum_{x\in\x_{I}}\frac{h_I^2(x){m_{21}}(x)}{x-y}=1+{m_{22}}(y),\\ 
&\sum_{x\in\x_{II}}\frac{h_{II}^2(x){m_{12}}(x)}{x-y}=1+{m_{11}}(y),\quad& 
\sum_{x\in\x_{II}}\frac{h_{II}^2(x){m_{22}}(x)}{x-y}={m_{21}}(y),
\endalignedat
$$
which is equivalent to our DRHP. 

This analogy suggests that continuous RHPs can be obtained as limits of DRHPs when $\x$ ``approximate''
$\Sigma$, or, {\it vice versa}, DRHPs can be obtained as limits of RHPs when the contour $\Sigma$ split into increasingly small parts passing through the points of $\x$ (the last observation is due to P.~Deift).
In \S8 we will provide an explicit example of a continuous RHP and a DRHP such that the discrete problem converges to the continuous one in a certain limit. 

Proposition 4.3 in our special situation takes the following form, cf. \S5.   

\proclaim{Proposition 6.1} Let $m$ be a solution of the DRHP stated above. Then the matrix $K=L(1+L)^{-1}$ has the form {\rm (}with respect to the splitting $\x=\x_I\cup\x_{II}${\rm )}
$$
\multline
K(x,y)=\\ \left[\matrix
\frac{h_I(x)h_I(y)(-m_{11}(x)m_{21}(y)+m_{21}(x)m_{11}(y))}
{x-y}&
\frac{h_I(x)h_{II}(y)(m_{11}(x)m_{22}(y)-m_{21}(x)m_{12}(y))}
{x-y}\\
\frac{h_{II}(x)h_I(y)(m_{22}(x)m_{11}(y)-m_{12}(x)m_{21}(y))}
{x-y}&
\frac{h_{II}(x)h_{II}(y)(-m_{22}(x)m_{12}(y)
+m_{12}(x)m_{22}(y))}
{x-y}
\endmatrix\right]
\endmultline
$$
where the indeterminacies of type $\frac 00$ on the diagonal are removed by the L'Hospital rule:
$$
K(x,x)=\cases h_I^2(x)(-m'_{11}(x)m_{21}(x)+m'_{21}(x)m_{11}(x)), &x\in\x_I,\\
h_{II}^2(x)(-m'_{22}(x)m_{12}(x)+m'_{12}(x)m_{22}(x)), &x\in\x_{II}.
\endcases
$$
\endproclaim
\demo{Proof}
A direct application of Proposition 4.3. Note that $\det m(\ze)\equiv 1$ by Lemma 4.7, and
$$
(m^t)^{-1}=\left(\matrix m_{22}& -m_{21}\\ -m_{12}& m_{11}\endmatrix\right).
\qed
$$
\enddemo

\head 7. The discrete Bessel kernel
\endhead

Now we return to the problem stated in \S2. As was explained in the previous section, the kernel $K=L(1+L)^{-1}$, where $L$ is given by (2.1), can be expressed via the solution of DRHP of the special form discussed above with 
$$
\gathered
\x=\Z'=\Z+\frac 12,\quad \x_I=\Z'_+=\left\{\frac 12,\frac 32,\dots\right\},\quad \x_{II}=\Z'_-=\left\{-\frac12,-\frac 32,\dots\right\},\\
h_I(x)=\frac{\theta^{\frac x2}}{\Gamma(x+\frac 12)}\,,\quad h_{II}(x)=\frac{\theta^{-\frac x2}}{\Gamma(-x+\frac 12)}\,.
\endgathered
$$

It will be more convenient for us to use the parameter $\eta=\sqrt{\theta}$ instead of $\theta$; then 
$$
h_I(x)=\frac{\eta^{x}}{\Gamma(x+\frac 12)}\,,\quad h_{II}(x)=\frac{\eta^{-x}}{\Gamma(-x+\frac 12)}\,.
\tag 7.1
$$

One way to extract the information about a solution of a RHP is to reduce the problem to a RHP with a constant (or not depending on a certain parameter) jump matrix. Then one can compare the solution with its derivatives with respect to the complex variable or with respect to a parameter, see, e.g., \cite{KBI, Ch. XV}. We will employ this approach for our DRHP. 

The asymptotic condition of the DRHP requires that $m(u)\sim 1$ as $u\to\infty$. Assume now that the next terms of the  asymptotic expansion of the solution $m(u)$ at infinity are given by the relation
$$
m(u)= \left(\matrix 1&0\\0&1\endmatrix\right)+\frac 1u
\left(\matrix \alpha&\beta\\ \gamma&\delta\endmatrix\right)+O\left(\frac 1{u^2}\right).
$$ 
Here $\al,\be,\ga,\de$ are some functions of $\eta$. The asymptotics here and below is understood in the sense that $|u|\to\infty$ so that $\operatorname{dist}(u,\x)$ is bounded from zero.

Note that the DRHP has an obvious symmetry: since $\x_I=-\x_{II}$ and $h_I(x)=h_{II}(-x)$, the changes 
$$
u\leftrightarrow -u, \quad \left(
\matrix m_{11}(u) & m_{12}(u)\\ m_{21}(u) &m_{22}(u)\endmatrix\right)\longleftrightarrow
\left(
\matrix m_{22}(-u) &- m_{21}(-u)\\- m_{12}(-u)& m_{11}(-u)\endmatrix\right)\,
\tag 7.2
$$
do not affect the problem (the minus signs in the off--diagonal blocks appeared because of the relation
$\res f(u)=-\operatorname{Res}\limits_{u=-x}f(-u)$). This
immediately implies that $\alpha=-\delta$, $\beta=\gamma$, so that
$$
m(u)= \left(\matrix 1&0\\0&1\endmatrix\right)+\frac 1u
\left(\matrix \alpha&\beta\\ \beta&-\alpha\endmatrix\right)+O\left(\frac 1{u^2}\right).
\tag 7.3
$$ 
Consider a new matrix 
$$
n(u)=m(u)\left(\matrix \eta^u&0\\0&\eta^{-u}\endmatrix\right).
$$
Observe that $n(u)$ is invertible for $u\in\C\setminus \x$, because $m(u)$ is invertible by Lemma 4.7. 

If $m(u)$ is the solution of the DRHP with $h$'s given by (7.1) and asymptotic behavior given by (7.3) then
$n(u)$ is the solution of the DRHP with
$$
\wt h_I(x)=\frac{1}{\Gamma(x+\frac 12)}\,,\quad \wt h_{II}(x)=\frac{1}{\Gamma(-x+\frac 12)}\,,
\tag 7.4
$$ 
and asymptotics
$$
n(u)= \left(\left(\matrix 1&0\\0&1\endmatrix\right)+\frac 1u
\left(\matrix \alpha&\beta\\ \beta&-\alpha\endmatrix\right)+O\left(\frac 1{u^2}\right)\right)\left(\matrix \eta^u&0\\0&\eta^{-u}\endmatrix\right),\quad u\to\infty.
\tag 7.5
$$ 
Note now that the jump matrix data of the new problem given by (7.4) do not depend on $\eta$, so that the derivative $\frac {\partial n}{\partial \eta}=\partial_\eta n$ satisfies the same DRHP with possibly different asymptotics. By Lemma 4.5, the matrix $(\partial_\eta n) n^{-1}$ has no singularities in $\C$.

The asymptotics of $(\partial_\eta n)n^{-1}$ can be easily computed from (7.5):
$$
\gathered
\partial_\eta n(u)= \frac u\eta  \left(\left(\matrix 1&0\\0&1\endmatrix\right)+\frac 1u
\left(\matrix \alpha&\beta\\ \beta&-\alpha\endmatrix\right)+O\left(\frac 1{u^2}\right)\right)\left(\matrix 1&0\\0&-1\endmatrix\right)\left(\matrix \eta^u&0\\0&\eta^{-u}\endmatrix\right),\\
n^{-1}(u)= \left(\matrix \eta^{-u}&0\\0&\eta^{u}\endmatrix\right)\left(\left(\matrix 1&0\\0&1\endmatrix\right)-\frac 1u
\left(\matrix \alpha&\beta\\ \beta&-\alpha\endmatrix\right)+O\left(\frac 1{u^2}\right)\right),\\
(\partial_\eta n)n^{-1}(u)=\frac 1\eta\left(\matrix u&-2\be\\2\be
&-u\endmatrix\right)+O\left(\frac 1u\right).
\endgathered
$$
Note that $\alpha$ has disappeared from the asymptotics. By the Liouville theorem we conclude that
$$
\partial_\eta n(u)\equiv \frac 1\eta\left(\matrix u&-2\be\\2\be
&-u\endmatrix\right) n(u).
\tag 7.6
$$
This relation already provides certain information about $n$. Indeed, it easily implies that the matrix elements of $n$ satisfy second order linear differential equations in $\eta$. However, the coefficients of this equation are still unknown
--- they are expressed in terms of the function $\beta(\eta)$.

The second trick which is commonly used in conventional RHPs in such a situation is differentiation with respect to the variable, if the jump matrix does not depend on the variable. Unfortunately, we cannot reduce our DRHP to one with jump matrix not depending on $u$. However, this obstacle can be overcome in the following way. 

Introduce a new matrix 
$$
p(u)=n(u)\left(\matrix \frac 1{\Gamma(u+\frac 12)}&0\\ 0 &\frac 1{\Gamma(-u+\frac 12)}\endmatrix\right).
$$
If $m(u)$ satisfies the original DRHP, it can be easily seen that $p(u)$ is holomorphic in $\C$. The residue condition of the DRHP and the asymptotics (7.5) in terms of $p(u)$ take the form
$$
p(x)=(-1)^{x-\frac 12}p(x)\left(\matrix 0&1\\1&0\endmatrix\right),\quad x\in\x,
\tag 7.7
$$
$$
p(u)= \left(\left(\matrix 1&0\\0&1\endmatrix\right)+\frac 1u
\left(\matrix \alpha&\beta\\ \beta&-\alpha\endmatrix\right)+O\left(\frac 1{u^2}\right)\right)\left(\matrix \frac{\eta^u}{\Gamma(u+\frac 12)}&0\\0&\frac{\eta^{-u}}{\Gamma(-u+\frac 12)}\endmatrix\right)
\tag 7.8
$$
as $u\to\infty$.  

Note that the relation (7.7) implies that the matrix $p(x)$ is {\it degenerate} at the points $x\in\x$.

Consider the matrix 
$$
\wt p(u)=\left(\matrix \wt p_{11}(u)&\wt p_{12}(u)\\ \wt p_{21}(u) &\wt p_{22}(u)\endmatrix\right)=
\left(\matrix  p_{11}(u+1)&- p_{12}(u+1)\\ - p_{21}(u-1) & p_{22}(u-1)\endmatrix\right).
$$
Clearly, $\wt p$ satisfies the condition (7.7), and (7.8) implies that, as $u\to\infty$,
$$
\gathered
\wt p_{11}(u)=(1+\alpha u^{-1}+O(u^{-2}))\,\frac{\eta^{u+1}}{\Gamma(u+\frac 32)}=\frac{\eta^u}{\Gamma(u+\frac 12)}\,O(u^{-1}),\\
\wt p_{12}(u)=-(\beta u^{-1}+O(u^{-2}))\,\frac{\eta^{-u-1}}{\Gamma(-u-\frac 12)}=\frac{\eta^{-u}}{\Gamma(-u+\frac 12)}\,\left(\frac{\be}\eta +O(u^{-1})\right),\\
\wt p_{21}(u)=-(\beta u^{-1}+O(u^{-2}))\,\frac{\eta^{u-1}}{\Gamma(u-\frac 12)}=\frac{\eta^u}{\Gamma(u+\frac 12)}\,\left(-\frac\be\eta +O(u^{-1})\right),\\
\wt p_{22}(u)=(1-\alpha u^{-1}+O(u^{-2}))\,\frac{\eta^{-u+1}}{\Gamma(-u+\frac 32)}=\frac{\eta^{-u}}{\Gamma(-u+\frac 12)}\,O(u^{-1}).
\endgathered
$$
This means that the matrix
$$
\wt n(u)=\wt p(u)\left(\matrix {\Gamma(u+\frac 12)}&0\\ 0 &{\Gamma(-u+\frac 12)}\endmatrix\right)
$$
satisfies the same DRHP as $n(u)$, and asymptotically
$$
\wt n(u)n^{-1}(u)\sim \left(\matrix 0&\frac \be\eta\\ -\frac\be\eta&0\endmatrix\right), \quad u\to\infty.
$$
By the Liouville theorem we get
$$
\wt n(u)= \left(\matrix 0&\frac \be\eta\\ -\frac\be\eta&0\endmatrix\right)n(u).
$$
In the language of $p$ and $\wt p$ this implies that
$$
\gathered
\wt p_{11}(u)=p_{11}(u+1)=\frac\be\eta\, p_{21}(u),\quad
\wt p_{21}(u)=-p_{21}(u-1)=-\frac\be\eta\, p_{11}(u),\\
\wt p_{12}(u)=-p_{12}(u+1)=\frac\be\eta\, p_{22}(u),\quad
\wt p_{22}(u)=p_{22}(u-1)=-\frac\be\eta\, p_{12}(u),
\endgathered
\tag 7.9
$$
for all $u\in\C$. Hence, $\left(\frac \be\eta\right)^2=1$, and $\be=\pm \eta$. Then the equation (7.6) takes the form
$$
\partial_\eta n(u)=\left(\matrix \frac u\eta &\mp 2\\\pm 2
&-\frac u\eta\endmatrix\right) n(u).
\tag 7.10
$$
Note that the change of signs of the off--diagonal elements of $n$ swaps the solutions of (7.10) corresponding to $\beta=\eta$ and $\beta=-\eta$, so it is enough to consider the case $\beta=-\eta$. Then (7.10) leads to the following relations
$$
\partial_\eta n_{11}(u)=\frac u\eta n_{11}(u)+2n_{21}(u),\quad \partial_\eta n_{21}(u)=-\frac u\eta n_{21}(u)-2n_{11}(u).
\tag 7.11
$$
This gives second order differential equations on
$n_{11}(u)$, $n_{12}(u)$:
$$
\left(\partial^2_\eta-\frac{u(u-1)}{\eta^2}+4\right)n_{11}(u)=0,\quad \left(\partial^2_\eta-\frac{u(u+1)}{\eta^2}+4\right)n_{21}(u)=0.
$$
General solutions of these equations are expressed through the Bessel functions
$$
\gathered
n_{11}(u)=\sqrt{\eta}\left(\const_1J_{u-\frac 12}(2\eta)+\const_2J_{-u+\frac 12}(2\eta)\right),\\ 
n_{21}(u)=\sqrt{\eta}\left(\const_3J_{u+\frac 12}(2\eta)+\const_4J_{-u-\frac 12}(2\eta)\right).
\endgathered
$$
The first line of (7.9) implies that (recall that $\beta=-\eta$)
$$
\const_1=-\const_3,\quad \const_2=-\const_4.
$$
The well--known differentiation formulas for the Bessel functions read
$$
\aligned
\partial_\eta(\sqrt{\eta}\,J_{\nu}(2\eta))=& \frac{\nu+\frac12}{\sqrt{\eta}} \,J_{\nu}(2\eta)-2\sqrt{\eta}\,J_{\nu+1}(2\eta)\\
=&\frac{-\nu+\frac12}{\sqrt{\eta}} \,J_{\nu}(2\eta)+2\sqrt{\eta}\,J_{\nu-1}(2\eta).
\endaligned
$$
Using these formulas in the first equation of (7.11) we see that it is equivalent to 
$$
4\const_2 \sqrt{\eta}\,J_{-u-\frac 12}(2\eta)=0.
$$
Hence, $\const_2=0$ and
$$
n_{11}(u)=\const\sqrt{\eta}\,J_{u-\frac 12}(2\eta), \quad n_{21}(u)=-\const\sqrt{\eta}\,J_{u+\frac 12}(2\eta)
$$
with the same constant (depending on $u$). Since $J_\nu(2\eta)\sim \eta^\nu/\Gamma(\nu+1)$ as $\nu\to\infty$, the asymptotic relation (7.5) implies that as $u\to\infty$,
$$
n_{11}(u)\sim\Gamma\left(u+\frac 12\right)\sqrt{\eta}\,J_{u-\frac 12}(2\eta), \quad n_{21}(u)\sim-\Gamma\left(u+\frac 12\right)\sqrt{\eta}\,J_{u+\frac 12}(2\eta).
$$
Moreover, $n_{11}(u)$ and $n_{21}(u)$ have the same singularities  as the corresponding expressions in the right--hand sides above (simple poles at the points of $\x_{II}$). Hence,
$$
n_{11}(u)=\Gamma\left(u+\frac 12\right)\sqrt{\eta}\,J_{u-\frac 12}(2\eta), \quad n_{21}(u)=-\Gamma\left(u+\frac 12\right)\sqrt{\eta}\,J_{u+\frac 12}(2\eta).
$$
Similar manipulations with $n_{12}$ and $n_{22}$ show that
$$
n_{12}(u)=\Gamma\left(-u+\frac 12\right)\sqrt{\eta}\,J_{-u+\frac 12}(2\eta), \quad n_{22}(u)=\Gamma\left(-u+\frac 12\right)\sqrt{\eta}\,J_{-u-\frac 12}(2\eta).
$$
Thus,
$$
m(u)=\sqrt{\eta}\left(\matrix J_{u-\frac 12}(2\eta)&J_{-u+\frac 12}(2\eta)\\
-J_{u+\frac 12}(2\eta)&J_{-u-\frac 12}(2\eta)\endmatrix\right)\left(\matrix {\eta^{-u}\Gamma(u+\frac 12)}&0\\ 0 &{\eta^{u}\Gamma(-u+\frac 12)}\endmatrix\right).
$$
It is easily verified that this matrix satisfies the original DRHP (the verification uses the symmetry relation
$J_{-n}=(-1)^nJ_n$, $n\in \Z$, for the Bessel functions).
Now it is immediately seen that the kernel $K(x,y)$ of Proposition 6.1 constructed from the matrix $m$ above coincides with the discrete Bessel kernel (2.2), and the Theorem 2.1 is proved.
\qed 

It is worth noting what happens if we choose $\be=\eta$ above. Then we get    
$$
\wh m(u)=\sqrt{\eta}\left(\matrix J_{u-\frac 12}(2\eta)&-J_{-u+\frac 12}(2\eta)\\
J_{u+\frac 12}(2\eta)&J_{-u-\frac 12}(2\eta)\endmatrix\right)\left(\matrix {\eta^{-u}\Gamma(u+\frac 12)}&0\\ 0 &{\eta^{u}\Gamma(-u+\frac 12)}\endmatrix\right),
$$
and this matrix does not satisfy the DRHP. However, it satisfies the problem up to the change of sign in the residue condition. This change of sign is equivalent to the multiplication of $h_I$ and $h_{II}$ by $\sqrt{-1}$, or to the change of sign of $L$. Thus, the kernel $\wh K$ constructed from the matrix $\wh m$ is connected with $L$ by the relation $$
\wh K=\frac L{L-1}.
$$

\head 8. Z--measures
\endhead
It seems that the DRHP considered in \S7 does not admit a scaling limit transition to a continuous RHP. The purpose of this section is to provide another DRHP for which such a limit exists. 

This DRHP and its scaling limit describe the so--called {\it z--measures on partitions} and the spectral decomposition of the {\it generalized regular representations} of the infinite symmetric group. We refer the reader to the papers \cite{KOV}, \cite{BOl1}, \cite{BOl2}, \cite{BOl3}, where this material is thoroughly explained.

We start with the description of the DRHP. As in \S7, we take $\x=\Z'$, $\x_I=\Z'_+$, $\x_{II}=\Z'_-$. The functions $h_I$, $h_{II}$ are as follows (we add the superscript `d' to the notation for the discrete problem and the superscript `c' to the notation for the continuous problem) 
$$
\gathered
h_I^d(x)=\frac{(zz')^\frac 14\xi^{\frac x 2}(1-\xi)^{\frac{z+z'}2}\sqrt{(z+1)_{x-\frac 12}(z'+1)_{x-\frac 12}}}{\Gamma(x+\frac 12)}\,,\\
h_{II}^d(x)=\frac{(zz')^\frac 14\xi^{-\frac x 2}(1-\xi)^{-\frac{z+z'}2}\sqrt{(-z+1)_{-x-\frac 12}(-z'+1)_{-x-\frac 12}}}{\Gamma(-x+\frac 12)}\,.
\endgathered
$$
Here $z$ and $z'$ are two complex parameters such that 
$(z+k)(z'+k)>0$ for all $k\in\Z$; $\xi\in(0,1)$ is also a parameter; $(a)_k=\Gamma(a+k)/\Gamma(a)$ is the Pochhammer symbol. Note that the limit 
$$
z,z'\to\infty,\quad\xi\to +0,\quad \xi\sqrt{zz'}=\theta
$$ 
brings us to the DRHP considered in \S7.

The continuous RHP in the notation of \S5 has the form
$$
\Sigma=\R\setminus\{0\},\quad\Sigma_I=\R_+,\quad\Sigma_{II}=\R_-,
$$
$$
\gathered
h_I^c(x)=\frac{(zz')^{\frac 14}}{\sqrt{\Gamma(z+1)\Gamma(z'+1)}}\,x^{\frac{z+z'}2}e^{-\frac x2},\\
h_{II}^c(x)=\frac{(zz')^{\frac 14}}{\sqrt{\Gamma(-z+1)\Gamma(-z'+1)}}\,(-x)^{-\frac{z+z'}2}e^{\frac x2}.
\endgathered
$$ 
It is not difficult to see that
$$
\lim_{\xi\to 1}h_I^d([(1-\xi)x])=h_I^c(x),\quad
\lim_{\xi\to 1}h_{II}^d([(1-\xi)x])=h_{II}^c(x).
$$
Here $[y]$ denotes the integer which is closest to $y$. 

In such a situation it is natural to say that the DRHP {\it approximates} the continuous RHP as $\xi\to 1$.

The solution for the DRHP was obtained in \cite{BOl2}:
$$
m^d(u)=\left(\matrix F\left(-z,-z';u+\tfrac12;\tfrac \xi{\xi-1}\right)&
\frac {\sqrt{zz'\xi}}{1-\xi}\,\frac{F\left(1+z,1+z';-u+\tfrac 32;\tfrac  \xi{\xi-1}\right)}{-u+\tfrac 12}\\
-\frac {\sqrt{zz'\xi}}{1-\xi}\,\frac{F\left(1-z,1-z';u+\tfrac 32;\tfrac  \xi{\xi-1}\right)}{u+\tfrac 12}&F\left(z,z';-u+\tfrac12;\tfrac \xi{\xi-1}\right)\endmatrix\right).
$$ 
Here $F(a,b;c;v)$ is the Gauss hypergeometric function. 

The solution for the continuous RHP was computed in different language in \cite{B1}, \cite{B2}, see also \cite{BOl1}:
$$
\multline
m^c(u)=\\ \left(\matrix u^{-\frac{z+z'+1}2}e^{\frac u2}W_{\frac{z+z'+1}2,\frac{z-z'}2}(u)&
\sqrt{zz'}\,(-u)^{\frac{z+z'-1}2}e^{-\frac u2}W_{\frac{-z-z'-1}2,\frac{z-z'}2}(-u)\\
-\sqrt{zz'}\,u^{-\frac{z+z'+1}2}e^{\frac u2}W_{\frac{z+z'-1}2,\frac{z-z'}2}(u)&(-u)^{\frac{z+z'-1}2}e^{-\frac u2}W_{\frac{-z-z'+1}2,\frac{z-z'}2}(-u)\endmatrix\right)
\endmultline
$$
where $W_{\kappa,\mu}(v)$ is the Whittaker function.

The correlation kernels $K$ corresponding to $m^d(u)$ and $m^c(u)$ are called the {\it hypergeometric kernel} and the {\it Whittaker kernel}, respectively. 

The convergence of the DRHP to the continuous RHP is established immediately using the relation
$$
\lim_{x\to+\infty}F\left(a,b;x;1-\frac xy \right)
=y^{\frac{a+b-1}2}e^{\frac y2}W_{\frac{-a-b+1}2,\frac{a-b}2}(y).
$$
We have
$$
\lim_{\xi\to 1} m^d((1-\xi)u)=m^c(u),
$$
see \cite{BOl2}, \cite{BOl3} for an explanation of the  representation theoretic meaning of this limit transition. 

It is worth noting that for the discrete problem described in this section the relation $K=L(1+L)^{-1}$ holds for all values of parameters, and $L$ and $K$ are bounded operators in $\ell^2(\Z')$. However, in the continuous case the operator $L$ becomes unbounded if $|\Re (z+z')|\ge 1$, and one should be careful to define $(1+L)^{-1}$. The spectral analysis of the kernels $L$ and $K$ in the continuous case has been done in \cite{Ol}, see also \cite{BOl1}.

\bigskip

\noindent Department of Mathematics, The University of
Pennsylvania, Philadelphia, PA 19104-6395, U.S.A. \newline and 
\newline Dobrushin Mathematics Laboratory, Institute for Problems of Information Transmission, Bolshoy Karetny 19, Moscow 101447, Russia. 

\noindent E-mail address:
{\tt borodine\@math.upenn.edu}

\enddocument